\documentclass[12pt]{article}\synctex=1 
\usepackage{amsfonts}
\usepackage{amssymb}
\usepackage{amsmath}
\usepackage{amsthm}
\usepackage{geometry}
\usepackage{commath}
\usepackage{booktabs}
\usepackage{lscape}
\usepackage{xr}
\usepackage[authoryear]{natbib}

\newcommand{\myref}[2]{\hyperref[#1]{#2}}
\numberwithin{equation}{section}

\usepackage{latexsym}
\usepackage{graphicx} 
\usepackage{enumerate}
\usepackage{setspace}
\usepackage[toc,title,titletoc,header]{appendix}
\usepackage[bottom]{footmisc}   
\usepackage[pdfborder={0 0 0}]{hyperref}

\usepackage{lineno}
\setpagewiselinenumbers

\newtheorem{theorem}{Theorem}[section]

\theoremstyle{definition}

\theoremstyle{remark}

\newcounter{assumptionM}
\newcounter{assumptionA}
\def\theassumptionM{M.\arabic{assumptionM}}
\def\theassumptionA{A.\arabic{assumptionA}}
\newcommand{\E}[1]{\mathbb{E}\left[#1\right]}

\usepackage[flushleft]{threeparttable}
\hypersetup{
  colorlinks=true,linkcolor=blue, citecolor=blue, filecolor=blue,
  linkbordercolor={0 1 1}, citebordercolor={1 0 0},
  pdftitle={A relaxed approach to estimating large portfolios},
  pdfauthor={Caner, Ulasan, Callot, Onder},
  pdfcreator={\LaTeX\ with package \flqq hyperref\frqq},
  pdfsubject={},
}

\begin{document}

\sloppy

\title{A Nodewise Regression Approach to Estimating Large Portfolios}

\author{Laurent Callot
\footnote{Amazon Research, Berlin, email: lcallot@amazon.com.}\and Mehmet Caner\footnote{
Department of Economics, North Carolina State University, email: mcaner@ncsu.edu}\and A.\"{O}zlem \"{O}nder
\footnote{Department of Economics, Ege University, email: ozlem.onder@ege.edu.tr.}\and Esra Ulasan
\footnote{Department of Economics, Ege University, email: esra.ulasan@ege.edu.tr}\thanks{We thank NYU-Stern Statistics, SOFIE 2017. We also thank editor Todd Clark, an associate editor, and two anonymous referees for comments that substantially improved the paper.
We thank Professor Anna O. Soter and Nancy House for professional text editing.}}

\date{\today}
\maketitle

\begin{abstract}

{\small
This paper investigates the large sample properties of the variance, weights, and risk of high-dimensional portfolios where the inverse of the covariance matrix of excess asset returns is estimated using a technique called nodewise regression. 
Nodewise regression provides a direct estimator for the inverse covariance matrix using the Least Absolute Shrinkage and Selection Operator (Lasso) of \cite{tibshirani94} to estimate the entries of a sparse precision matrix.
We show that the variance, weights, and risk of the global minimum variance portfolios and the Markowitz mean-variance portfolios are consistently estimated  with more assets than observations.
We show, empirically, that the nodewise regression-based approach performs well in comparison to factor models and shrinkage methods.}


\end{abstract}
\newpage

\section{Introduction}

Accurately estimating the variance, weights, and risk is crucial to the formation of portfolios of financial assets. These three quantities are functions of the inverse covariance matrix, or precision matrix, of excess returns. In this paper, we establish the consistency of portfolio variances, weights, and risks when the precision matrix is estimated using the nodewise-regression method of \cite{meinshausen06}, even where the number of assets is greater than the number of observations. We define the number of assets as $p$, and the sample size or the time span of the portfolio as $n$.

Nodewise regression provides a direct estimator of the $p \times p$ sparse precision matrix by using only $p$ Lasso regressions. The entries of the precision matrix are given by an exact formula as the solution to a regression problem. The inverse covariance matrix is therefore directly estimated and not given as the inverse of some estimated covariance matrix. The entries of the precision matrix are sparse in high-dimensions. This does not imply that the corresponding covariance matrix is sparse. \cite{meinshausen06} show that the sparsity pattern of the precision matrix can be tied to the optimal predictor for the excess asset returns. We discuss  this point further in section 2.3 below. \cite{vandegeer14}, and \cite{caner14} use nodewise regression to build asymptotically honest confidence intervals for high-dimensional parameter vectors. In a time series context, \cite{chang18} build confidence intervals around each entry in a precision matrix estimated by nodewise regression.

There are two main alternative strategies to sparsity for estimating high-dimensional precision matrices in a portfolio optimization context, namely,  factor models and shrinkage. In the following paragraphs, we discuss papers from classes of models that are most closely related to ours, focusing on papers establishing theoretical properties for the portfolios formed using the estimator they propose.

Factor models assume that the excess asset returns are driven by a small number of latent or observed factors. The performance of these models is documented in \cite{fan2011}, and \cite{fan2013large}. In this class of models, the factor-based covariance matrix of excess returns is estimated. The factor-based covariance matrix can be inverted, due to the low-dimensional factor structure. \cite{fan2016}, and \cite{ait2016} investigate the accuracy of precision matrices estimated by using factor models, but do not analyze portfolio variances and weights when $p>n$.

\cite{fan08}, propose a covariance matrix estimator assuming observable factors and a diagonal residual covariance matrix, that is, assuming conditional sparsity. They establish the rates of convergence of the portfolio risk in large portfolios using their covariance matrix estimator.
\cite{fan15}, propose the principal orthogonal complement thresholding (POET) estimator where the factors are unobserved. In a high-dimensional setting, \cite{li17} provide valuable findings on the optimal global variance ratio via factor models, and \cite{ao17} on portfolio selection using the Lasso when asset returns are heteroskedastic.

An alternative strategy for estimating precision matrices in large dimensions, is based  on regularizing the sample covariance matrix of asset returns and then inverting it. \cite{lediot03,ledoit04} propose shrinkage approaches using convex combinations of the sample covariance matrix of excess asset returns, and  an identity matrix. \cite{ledoit12} propose a non-linear shrinkage of the sample eigenvalues for the case where $p<n$, and \cite{ledoit2015spectrum} extend  this approach to the case where $ p>n $. Because \cite{ledoit04, ledoit2015spectrum} are primarily concerned with covariance estimation, they do not study the theoretical properties of portfolios formed using their estimators of the precision matrix.

This paper makes the following contributions to the literature on high-dimensional portfolios. First, we establish consistency results for the variance of the global minimum variance portfolio and the \cite{markowitz52} portfolio when the number of assets is larger than the sample size.
To the best of our knowledge, there are no equivalent theoretical results for factor models or shrinkage methods.
We prove that the weights of the global minimum variance and \cite{markowitz52} portfolios can be consistently estimated when the number of the assets is larger than the sample size. 
These results are established for both fixed and growing exposure, where exposure is defined as the $l_1$ norm of the portfolio weights. 
We are only aware of similar results in the context of factor models with constant exposure \citep{fan15,fan08}, in which case the rate of convergence is identical to ours.

Finally, we show that the risk of global minimum variance and the \cite{markowitz52} portfolios are consistently estimated when the number of assets is larger than the sample size under both constant and growing exposure. To the best of our knowledge, similar results have only been established in the constant exposure case with factors models \citep{fan15,fan08}, and have not been established using shrinkage methods.

The rest of this paper is organized as follows. In Section \ref{sec:Nodewise}, we introduce the nodewise regression with $\ell_1$ penalty as well as the approximate inverse of the empirical Gram matrix and its asymptotic properties. In Section \ref{main}, we establish our main theoretical results on the convergence rates of the estimator of the variance, the weights, and the risk of global minimum portfolios and mean-variance portfolios based on the nodewise regression estimator of the precision matrix. In Section \ref{sec:simulation}, we report simulation results comparing nodewise-based portfolios to portfolios based on alternative estimators of the inverse covariance matrix. In Section \ref{sec:per}, we report the results of an empirical application, and in section \ref{sec:conclusion}, we present our  conclusions. All the proofs are provided in the Appendix.

\section{The Lasso for Nodewise Regression}\label{sec:Nodewise}

In this section, we  introduce the nodewise regression estimator developed by \cite{meinshausen06} and \cite{vandegeer14}. We  begin by deriving the exact formula for the precision matrix, then discuss its estimation by nodewise regression, and finally discuss the advantages of the nodewise regression approach and contrast it with alternative methods. Before proceeding, we introduce the following notation. 

Throughout the paper, $\| \nu \|_{\infty} \|,  \nu \|_1, \| \nu \|_2$ denote the sup,  $l_1$, and Euclidean norm of a generic vector $\nu$, respectively. $\|A\|_{\infty}$ is the sup norm for matrices, and denotes the largest absolute value for a generic matrix  $A$, $A \ge 0$, implies a positive semidefinite matrix.

Define $r_t$ as  a $p \times 1$ vector of excess asset returns at time $t=1,\cdots, n$ and $\E{r_t} = \mu :=(\mu_1, \cdots, \mu_j, \cdots, \mu_p)'$. Assume $\| \mu \|_{\infty} \le C < \infty$, where $C$ is a positive constant. 
Let $ \Sigma\in\mathbb{R}^{p\times p}$ be the full rank covariance matrix of excess asset returns where $ \Sigma :=\mathbb{E}\left[(r_{t}-\mu)(r_{t}-\mu)^{\prime}\right]$.
The sample covariance matrix is $\hat{\Sigma} := n^{-1} \sum_{t=1}^n (r_t - \bar{r}) (r_t - \bar{r})'$, with $\bar{r}:= n^{-1} \sum_{t=1}^n r_t $. 
For each $j=1, \cdots , p$, define the demeaned excess returns of asset $j$ as $r^*_j  := r_j - \bar{r}$, and let $\mathbf{r^*} :=\left[ r^*_{1},\ldots,r^*_{p}\right] $ be the matrix of excess asset returns of dimensions $n\times p$. $r^*_{j}$ is a column vector of length $n$, and denotes the $j^{th}$ column in $\mathbf{r^*}$. $\mathbf{r^*}_{-j}$ denotes all columns of $\mathbf{r^*}$ except for the $j^{th}$ one. Finally, define $\Theta := \Sigma^{-1}$.

\subsection{Derivation of the Precision Matrix Formula}

Let $\Sigma_{-j,-j}$ represent the $p-1 \times p-1$ sub-matrix of $\Sigma$ where the $j$th row and column have been removed.
Similarly, let $\Sigma_{j,-j}$  denote the $j$ th row of $\Sigma$ with its $j$ th element removed, and let $\Sigma_{-j,j}$ denote the $j$ th column of $\Sigma$ with its $j$ th element removed. 
Using the inverse formula for block matrices, the $j$ th main diagonal term is 
\begin{equation}
\Theta_{j,j} = ( \Sigma_{j,j} - \Sigma_{j,-j} \Sigma_{-j,-j}^{-1} \Sigma_{-j,j})^{-1},\label{pm0}
\end{equation}
and the $j$ th row of $\Theta$ with the $j$ th element removed is
\begin{equation}
\Theta_{j,-j} = - ( \Sigma_{j,j} - \Sigma_{j,-j} \Sigma_{-j,-j}^{-1} \Sigma_{-j,j})^{-1} \Sigma_{j,-j} \Sigma_{-j,-j}^{-1} = - \Theta_{j,j} \Sigma_{j,-j} \Sigma_{-j,-j}^{-1}.\label{pm1}
\end{equation}

We now show how  (\ref{pm0}) and (\ref{pm1}) are related to a linear regression. Let $r_{t,j}^*:= r_{t,j} - \frac{1}{n} \sum_{t=1}^n r_{t,j}$ denote the demeaned\footnote{\cite{vandegeer14} do not use demeaned data while \cite{chang18} do. Demeaning or not makes no noticeable difference in our empirical results.} excess returns of asset $j$ at time $t$, and let $r_{t,-j}^*$ be the vector of all demeaned returns except for the $j$ th one.

Define $\gamma_j$ as the value of the vector $\gamma$ of length $p-1$ that minimizes $\E{r_{t,j}^* - (r_{t,-j}^*)' \gamma}^2$ for all $t=1,\cdots,n$.  We get the solution 
\begin{equation}
\gamma_j = \Sigma_{-j,-j}^{-1} \Sigma_{-j,j},\label{pm2}
\end{equation} 
by using strict stationary of the data. Using the symmetry of $\Sigma$ and (\ref{pm2}), we can re-write   (\ref{pm1}) as 
\begin{equation}
\Theta_{j,-j} = - \Theta_{j,j} \gamma_j'.\label{pm3}
\end{equation}
Define $\eta_{t,j} := r_{t,j}^* - (r_{t,-j}^*)' \gamma_j$. By  (\ref{pm2}), we can verify that 
\begin{eqnarray}
\E{r_{t,-j}^* \eta_{t,j}} & = & \E{r_{t,-j}^* r_{t,j}^*} - \E{r_{t,-j}^* (r_{t,-j}^*)'} \gamma_j \nonumber \\
& = & \Sigma_{-j,j} - \Sigma_{-j,-j} \Sigma_{-j,-j}^{-1} \Sigma_{-j,j} = 0.\label{pm3a}
\end{eqnarray}

Gathering the results above, we can pose the following regression model with covariates orthogonal to errors, namely:
\begin{equation}
r_{t,j}^* = (r_{t,-j}^*)' \gamma_j + \eta_{t,j}.\label{pm4}
\end{equation}
From (\ref{pm3}) and (\ref{pm4}), we can see that $\Theta_{j,-j}$, and hence the row $\Theta_j$, is sparse if and only if $\gamma_j$ is sparse.

We now use the results above to derive a formula for $\Theta$. Using  (\ref{pm2}), (\ref{pm3a}), and (\ref{pm4}), we get 
\begin{eqnarray}
\Sigma_{j,j} = \E{r_{t,j}^*}^2 & = & \gamma_j' \Sigma_{-j,-j} \gamma_j + \E{\eta_{t,j}^2} \nonumber \\
& = & \Sigma_{j,-j} \Sigma_{-j,-j}^{-1} \Sigma_{-j,j} + \E{\eta_{t,j}}^2.\label{pm5}
\end{eqnarray}
Define $\tau_j^2 := \E{\eta_{t,j}^2}$. By (\ref{pm5}),
\begin{equation}
\tau_j^2 = \Sigma_{j,j} - \Sigma_{j,-j} \Sigma_{-j,-j}^{-1} \Sigma_{-j,j} = \frac{1}{\Theta_{jj}},\label{pm5a}
\end{equation}
where we use  (\ref{pm0}) for the second equality. Next, define a $p \times p$ matrix
\[ C := \left[ \begin{array}{cccc}
1 & -\gamma_{1,2} & \cdots & -\gamma_{1,p} \\
-\gamma_{2,1}& 1 & \cdots  & \cdots \\
\vdots & \vdots & \vdots & \vdots \\
-\gamma_{p,1} & - \gamma_{p,2}&  \cdots  & 1 
\end{array}
\right],\]
and define the diagonal matrix $T^{-2} := diag (\tau_1^{-2}, \cdots, \tau_p^{-2})$. Finally,   we get 
\begin{equation}
\Theta = T^{-2} C,\label{pm7}
\end{equation}
 since (\ref{pm5a}) establishes  that $\Theta_{j,j} = \frac{1}{\tau_j^2}$, and by (\ref{pm3}), $\Theta_{j,-j} = - \Theta_{j,j} \gamma_j' = \frac{-\gamma_j'}{\tau_j^2}$.

\subsection{Nodewise Regression Estimation of the Precision Matrix}\label{method}

This subsection borrows from \cite{vandegeer14}, and  \cite{chang18}  to introduce the nodewise regression estimation algorithm. For  each $j=1,\ldots,p$, the nodewise regression is defined as
\begin{equation}
\hat{\gamma}_{j}:=\underset{\gamma\in\mathbb{R}^{p-1}}{\text{argmin}} (\parallel r^*_{j}-\mathbf{r^*}_{-j}\gamma\parallel_{2}^{2}/n+2\lambda_{j}\parallel\gamma\parallel_{1}), \label{gamma}
\end{equation}
where $\hat{\gamma}_{j}=\left\lbrace \hat{\gamma}_{j,k}; k=1,\ldots,p, k\neq j\right\rbrace $ is the vector of length $ (p-1)$ of regression coefficient estimates which will be used to construct the estimate of the precision  matrix. $\lambda_{j}$ is a positive tuning parameter that  determines the size of the penalty on the parameters.  
Let $S_{j}:=\lbrace k;\: \gamma_{j,k}\neq0\rbrace$ be the set of non-zero estimates for row $\gamma_j$ in (\ref{pm3}), and let $s_{j}:=\vert S_{j} \vert$ be its cardinality.

The nodewise regression estimator $ \hat{\Theta}$ of the precision matrix ${\Theta}$ is constructed as follow. First, define
\[ \hat{C}:=
\begin{pmatrix}
1 & -\hat{\gamma}_{1,2} & \dots  & -\hat{\gamma}_{1,p} \\
-\hat{\gamma}_{2,1} & 1 & \dots  & -\hat{\gamma}_{2,p}  \\
\vdots & \vdots & \ddots & \vdots \\
-\hat{\gamma}_{p,1} & -\hat{\gamma}_{p,2}  &  \dots  & 1
\end{pmatrix}, \label{c}
\]
and write $\hat{T}^{2}:=diag(\hat{\tau}_{1}^{2},\ldots,\hat{\tau}_{p}^{2})$, the $p\times p$ diagonal matrix, with entries
\begin{equation}
\hat{\tau}_{j}^{2}:=\dfrac{\parallel r^*_{j}-\mathbf{r^*}_{-j}\hat{\gamma}_{j}\parallel^{2}_{2}}{n}+\lambda_{j}\parallel \hat{\gamma}_{j}\parallel_{1}.\label{tau}
\end{equation}
Then, define the approximate inverse $\hat\Theta:=\hat{T}^{-2}\hat{C}$ as the estimator of  (\ref{pm7}). Notice that while $\hat{\Sigma}$ is self-adjoint, $\hat{\Theta}$ is not. 

The penalty parameter $\lambda_j$ in (\ref{gamma}) is chosen by minimizing the Generalized Information Criterion  (GIC) of \cite{fan13}:
\begin{align}
  GIC(\lambda_j) := \log(\hat{\sigma}_{\lambda_j}^2) + | \hat{S}_j (\lambda_j)| \frac{\log(p)}{n}\log(\log(n)),\label{mbic}
\end{align}
where $\hat{\sigma}_{\lambda_j}^2:= \| r_j - {\bf r}_{-j} \hat{\gamma}_j  \|_2^2/n$ is the residual variance for asset $j$, and $|\hat{S}_j (\lambda_j)|$ represents the estimated number of nonzero parameters in the vector $\hat\gamma_j$. 
We use the GIC to select $\lambda_j$ as Corollary 1 and Theorem 2 of \cite{fan13} shows that the GIC selects the true model with probability approaching one both when $p>n$ and when $p \le n$.

The algorithm below summarizes the steps to compute $\hat{\Theta}$:
\begin{enumerate}
  \item Estimate $\hat{\gamma}_{j}$ for a given $\lambda_j$ by solving  (\ref{gamma}).
  \item Choose $\lambda_{j}$ using the GIC of (\ref{mbic}).
  \item Repeat steps 1-2 for $j=1,\cdots,p$.
  \item Compute $\hat{C}$ and $\hat{T}^{2}$. 
  \item Return the nodewise estimator of the precision matrix $\hat\Theta=\hat{T}^{-2}\hat{C}$. \label{algorithm}
\end{enumerate}

\subsection{Nodewise Regression and Optimal Prediction}

\cite{meinshausen06} establish a link between nodewise regression and the optimal linear prediction of excess asset returns under the assumption that the returns are normally distributed. Throughout this paper, we do not assume normality, but will do so only in this subsection.

Define a neighborhood $ne_j$ asset $j$ as being the smallest subset of $\{ 1, \cdots, p \} - \{ j \}$ such that $r_{t,j}^*$ is conditionally independent of all other assets outside the neighborhood.  Consider the optimal prediction of $r_{t,j}^*$ given a subset of $\{ r_{t,k}^*: k \in {\cal A}\}$, where ${\cal A} \subseteq \{1,2,\cdots, p\} - \{ j \}$.

\[ \gamma_j^* = argmin_{\gamma_j \gamma_{j,k} =0, k \notin {\cal A}} E [ r_{t,j}^* - \sum_{k \in \{1, 2, \cdots, p\}} \gamma_{j,k} r_{t,k}^*]^2.\]
 
This is  equation (2) in  \cite{meinshausen06}. Note that the set of non-zero coefficients of $\gamma_j^*$ defines the neighborhood $ne_j$ of $r_{t,j}^*$. This set is identical to the set of non-zero entries in the row entries in the precision matrix, thereby providing the link between nodewise regression and the optimal predictor of the returns. 
Similarly, \cite{yuan10}, ties the entries of the precision matrix to regression coefficients in the case of \textit{iid} random variables in a non-sparse setting as we have done in detail in Section 2.1.



\subsection{Comparison with Alternative Approaches}

Other techniques, such as factor-based methods, use a conditionally sparse covariance matrix and then invert it. There are two main underlying assumptions in a factor-model based precision matrix. First, a limited number of latent factors explain covariance between assets. Second, some correlations between shocks may be zero.  Factor model is a complementary approach since the dimension reduction is done in the covariance matrix, after which  this new structure is estimated and inverted.

Shrinkage-based methods also estimate covariance matrices before inverting them. The advantage of shrinkage-based methods is that they do not assume (conditional) sparsity of the covariance matrix. However their consistency properties have not been established, and \cite{ao18} show that they fail to achieve mean-variance efficiency.

There has been a growing interest in recent years in the estimation of covariance matrices and precision matrices in high dimensions. 
One approach is to use the projection method by \cite{fanj17}, and \cite{zhao14}. \cite{fanj17} provide a theory of estimation through using  the projection method, and \cite{zhao14} provide and estimation algorithm. Specifically, \cite{fanj17} show that under the Gaussian copula model, the sparsity of the precision matrix is equivalent to conditional independence of variable j from variable k,  given other variables. They introduce an estimator for $\Sigma$ based on rank estimation. The sample covariance matrix is replaced by a rank-based estimator which may not be positive semi-definite. They propose a method to make the sample covariance matrix positive-definite. 
Specifically,  the sample covariance is projected onto the cone of positive semi-definite matrices. This projection is inverted instead of the sample covariance. Large sample portfolio theory for the method of \cite{fanj17}, and \cite{ zhao14} has not yet been established. 

Other approaches include that of \cite{ek10}, where the properties of the Markowitz portfolio based on quadratic optimization, are established in a setting where $p$ is of the same magnitude as $n$, and with independent, normally distributed asset returns.

In a similar vein, \cite{li15} and \cite{yen16} consider quadratic portfolio optimization problems with linear constraints and $l_1, l_2$ penalties.  The portfolio weights are assumed to be  sparse, which is a more restrictive assumption than the assumption that the precision matrix is sparse as required for nodewise regression. Similarly, \cite{li15}, \cite{yen16}, also do not establish consistency or other large sample properties of the risk and the variance of portfolios based on their proposed method.
 
 \cite{dai17} analyze the asymptotic properties of the BARRA estimator which uses factor models to enable an analysis of large covariance matrix and precision matrix. There is, however, no connection to the financial issues that we analyze here.  \footnote{We thank two anonymous referees for pointing out these references. }
 
 In summary, our approach is based on a closed-form expression for each element in the sparse precision matrix, and we have been  able to establish the large sample properties of quantities that are key for portfolio formation even when $p>n$. Our method can be seen as a complementary approach to factor models and shrinkage methods. 
 In addition, with Gaussian data, it is possible to show that nodewise regression provides an optimal predictor of returns.

\section{Nodewise Regression Estimator in Large Portfolios}\label{main}

In this section, we establish the consistency of the estimators of the variance, risk, and weights of the global minimum variance and the \cite{markowitz52} portfolios. 

\subsection{Assumptions}\label{subsec:asym}
 To derive our theorems, we make the following assumptions.   Set $\bar{s} := \max_{ 1 \le j \le p} s_j$, where $s_j$ is the number of non-zero parameters in row $j$ introduced in Section 2.2.

{\bf Assumption 1}. {\it The $n \times p $ matrix of excess asset returns $\mathbf{r^*}$ has strictly stationary $\beta$ mixing rows with $\beta$ mixing coefficients satisfying
$\beta_k \le exp( -K_1 k^{\Xi_1})$ for any positive  $k$, with constants $K_1 > 0, \Xi_1>0$ that are independent of $p$ and $n$.}

{\bf Assumption 2}. {\it The smallest eigenvalue of $\Sigma$, $\Lambda_{min}$, is strictly positive and uniformly bounded away from zero. The maximum eigenvalue of $\Sigma$ is uniformly bounded away from infinity.}

{\bf Assumption 3}.{\it There exist constants $K_2> 0, K_3 > 1, 0 < \Xi_2 \le 2, 0 < \Xi_3 \le 2$ that are independent of $p$ and $n$ such that
\[ \max_{1 \le j \le p} E [ exp (K_2 | r^*_{tj} )|^{\Xi_2} ] \le K_3,\]
\[ \max_{1 \le j \le p} E [ exp (K_2 | \eta_{tj})|^{\Xi_3} ] \le K_3.\]}
{\bf Assumption 4}. {\it The following condition holds: $ \bar{s} \sqrt{\mathrm{log}p/n} = o(1).$}

Assumptions 1-3 are similar to the assumptions made in \cite{chang18}.  Assumption 1 allows for weak dependence in the data. As shown in \cite{chang18}, causal ARMA processes with continuous innovation distributions are $\beta$ mixing with exponentially decaying $\beta_k$. Certain stationary Markov chains also satisfy this condition. Stationary GARCH models with finite second moments and continuous innovation distributions  also satisfy Assumption 1. The details are provided in \cite{chang18}.

Assumption 2 requires the population covariance matrix to be nonsingular, assuming a strictly positive minimum eigenvalue rules out local to zero sequences. 

Assumption 3 restricts the tail behavior which allowed \cite{chang18} to use exponential tail inequalities to establish upper bounds on certain tail probabilities.

Assumptions 1 and 3 could be replaced with assumptions that allow for independent but non-identical data with bounded moments as discussed in Appendix B of \cite{caner14}.

Assumption 4 is a sparsity condition on the inverse of the covariance matrix of excess asset returns. It does not imply the sparsity of the covariance matrix of excess asset returns.  For example, if the excess asset returns have an autoregressive structure of order one, the covariance matrix will be non-sparse, but its inverse will be sparse.
Another example assumes that $\Sigma$ has a block diagonal or a Toeplitz structure $\Sigma_{i,j} = \rho^{|i-j|}$, $-1< \rho< 1$, with $\rho$ being the correlation among assets. In this case, the inverse of $\Sigma$ is sparse.

Regarding Assumption 4, note that when $p>n$, $\bar{s} =o(\sqrt{n/logp})$. In the simple case of $p=a n$, where $a>1$ is constant, $\bar{s}= o(\sqrt{n/log(an)})$  is growing with $n$ but is smaller than the maximum number of possible non-zeros $p-1$ in a row.

This sparsity assumption could be replaced with a weaker condition allowing for many small coefficients in the inverse without the need to impose sparsity when $p>n$ by defining $\Xi_r = \max_{1 \le j \le p} \sum_{k=1}^p |\gamma_{jk}|^r$ and assuming that $\lim_{ r \to 0} \Xi_r \le C < \infty$ when $r \to 0$ as in \cite{vandegeer16}. We do not pursue this extension in the present paper.

If the problem is such that $p<n$, then it is possible to assume that $\bar{s} = p-1$ (no sparsity in the inverse of the covariance matrix) in which case Assumption 4 can be written as $p \sqrt{logp/n} = o(1)$.    Note that our Assumption 4 has an extra $\sqrt{\bar{s}}$ factor compared to that of \cite{vandegeer14} which is due to  portfolio optimization problems.


\subsection{Optimal Portfolio Allocation and Risk Assessment}\label{subsec:port}


A portfolio is a set of weights $w=\left(w_{1},\ldots,w_{p}\right)^\prime\in\mathbb{R}^{p}$ representing the relative amount of wealth invested in each asset. A full investment constraint requires that weights should sum to $1$. It is written $w^\prime1_{p}=1$ where $ 1_{p}=(1, \ldots, 1)^\prime$. Throughout the paper, we assume short-selling is allowed and hence the value of individual weights can be negative.

\subsubsection{Global Minimum Variance Portfolio} \label{sgmp}

The global minimum variance portfolio is the set of weights $w_u$ that minimizes the portfolio's variance $w^\prime \Sigma w$. The portfolio weights are given by:
\[ \hat w_u=\text{arg} \underset{w}{\text{min}} (w^\prime \Sigma w), \quad \text{such that}  \quad w' 1_p = 1. \]

\noindent Define $A = 1_p' \Theta 1_p/p$, where $\Theta = \Sigma^{-1}$. $A$ is estimated by $\hat{A} = 1_p' \hat{\Theta} 1_p/p$. The global minimum variance is $\Phi_{G} = w_u' \Sigma w_u = (p A ) ^{-1}$ as shown in \cite{fan08} (equation (11)), and  is estimated by $\hat{\Phi}_{G} = (p \hat{A})^{-1}$.

The next result,  one of the main results of this paper, shows that the global minimum variance is consistently estimated by using the nodewise regression estimator of the inverse covariance matrix.
\begin{theorem}\label{thm1}
 {\it Under Assumptions 1-4, with $\lambda_j = O( \sqrt{logp}/n)$ uniformly in j,
\[   |\frac{\hat{\Phi}_{G}}{  \Phi_G} -1 | = O_p ( \bar{s} \sqrt{logp/n})= o_p(1).\]}
\end{theorem}

 To the best of our knowledge, establishing the consistency of estimate of a large portfolio's variance to its population counterpart without assuming factor structure, is a novel contribution. \cite{li17} derives a similar result in a factor model setting.

\subsubsection{Markowitz Mean-Variance Framework}\label{ssmkw}

The portfolio selection problem as defined by \cite{markowitz52}, is to find the portfolio with the smallest variance given a desired expected return $\rho_1$. At time $ t $, an investor determines the portfolio weights that minimize the mean-variance objective function:

\begin{equation}
  \hat w=\text{arg} \underset{w}{\text{min}} (w^\prime \Sigma w), \quad \text{subject to} \quad  w^\prime1_{p}=1 \quad  \mathrm{and} \quad  w^\prime\mu=\rho_1 \label{markowitz}.
\end{equation}

\noindent Define the terms $B = 1_p' \Sigma^{-1} \mu/p$ and $D= \mu' \Sigma^{-1} \mu/p$, and their estimates $\hat{B} =  1_p' \hat{\Theta} \hat{\mu}/p$ and
$\hat{D} = \hat{\mu}' \hat{\Theta} \hat{\mu}/p$, with $\hat{\mu} := \bar{r}= n^{-1} \sum_{t=1}^n  r_t$.

\noindent The optimal portfolio solution to the constrained optimization problem (\ref{markowitz}) is
\begin{equation}
w^{*'} \Sigma w^* = p^{-1} \left[ \frac{ A \rho_1^2 - 2 B \rho_1 + D}{A D - B^2}\right].
\label{consop}
\end{equation}

\noindent The estimate for the above optimal portfolio variance is
\[ p^{-1} \left[ \frac{ \hat{A} \rho_1^2 - 2 \hat{B} \rho_1 + \hat{D}}{\hat{A} \hat{D} - \hat{B}^2} \right].\]

Further, define the optimal portfolio variance as $\Psi_{OPV} := p^{-1} \left[\frac{ A \rho_1^2 - 2 B \rho_1 + D}{A D - B^2}\right]$ and its estimate $\hat{\Psi}_{OPV} := p^{-1} \left[\frac{ \hat{A} \rho_1^2 - 2 \hat{B} \rho_1 + \hat{D}}{\hat{A} \hat{D} - \hat{B}^2}\right]$.

The following Theorem shows that the optimal portfolio variance is consistently estimated using the nodewise regression estimator of the inverse covariance matrix.
\vspace{0.5cm}
\begin{theorem}\label{thm2}
{\it Under Assumptions 1-4, with $\lambda_j = O ( \sqrt{\mathrm{log}p/n})$ uniformly in j, with $ (AD-B^2)\ge C_1 > 0$ and $ (A \rho_1^2 - 2 B \rho_1 + D) \ge  C_1 $, where $C_1$ is a positive constant, and with $\rho_1$ uniformly bounded away from infinity, we get
\[ \left|\frac{\hat{\Psi}_{OPV}}{  \Psi_{OPV}} -1  \right| =  O_p ( \bar{s}  \sqrt{logp/n})=  o_p (1).\]
}
\end{theorem}

We show that the ratio of the estimated portfolio variance to its population counterpart ratio converges to one in probability, which is, to the best of our knowledge, a new result in the analysis of large portfolios. 

The restrictions $(A \rho_1^2 - 2 B \rho_1 + D) \ge  C_1 >0$ and $(AD-B^2)\ge C_1 > 0$ ensure that the variance of the optimal portfolio $\Psi_{OPV}$ in (\ref{consop}) is positive, and finite. 
Lemma A.4 provides a large sample analysis  of the terms $A, B$, and $D$.

\subsection{Portfolio Weights Estimation}\label{subsec:gross}

In this subsection, we establish the properties of the estimated weights of portfolios based on the nodewise regression estimator of the inverse covariance matrix.
The weights of the global minimum variance portfolio are given by:
\begin{equation}
w_u = \frac{\Sigma^{-1} 1_p}{1_p'  \Sigma^{-1} 1_p } = \frac{ \Sigma^{-1} 1_p/p}{A}, \label{gmw}
\end{equation}
and are estimated by:
\begin{equation}
\hat{w}_u = \frac{\hat{\Theta} 1_p }{1_p' \hat{\Theta} 1_p} =
\frac{\hat{\Theta} 1_p/p}{\hat{A}}. \label{gmt}
\end{equation}

Theorem \ref{thm3} shows that we can consistently estimate the weights of the global minimum variance portfolio in a high dimensional setting.

\begin{theorem}\label{thm3}
 {\it Under Assumptions 1-3 and the sparsity assumption $\bar{s}^{3/2} \sqrt{\mathrm{log}p/n}=o(1)$, with $\lambda_j = O (\sqrt{logp/n})$ uniformly in j, we have
\[ \|\hat{w}_u - w_u \|_1 = O_p ( \bar{s}^{3/2} \sqrt{\mathrm{log}p/n})=o_p(1).\]}
\end{theorem}

{\bf Remarks:}
\begin{enumerate}
\item The condition $\bar{s}^{3/2} \sqrt{\mathrm{log}p/n}=o(1)$ adds an extra $\sqrt{\bar{s}}$ to Assumption 4 to limit the $\ell_1$ norm of the approximation error for $\hat{\Theta}$. 

\item Consider the case where we assume  that  off-diagonal elements of $\hat{\Theta}$ are not sparse, that is, $\bar{s} = p-1$. In that case, $p$ has to be less than $n^{1/3}$, implying $p (logp)^{1/3} = o(n^{1/3})$ to satisfy the assumption
$(\bar{s})^{3/2} \sqrt{logp/n} = (p-1)^{3/2} \sqrt{logp/n} = o(1)$.

\item  $\|w_u \|_1$is the gross exposure of the portfolio. By Equation (A.52) (in the appendix) and a simple eigenvalue inequality, we would have $\|w_u \|_1 = O (\sqrt{\bar{s}})$: the gross exposure may grow with $n$. Theorem \ref{thm3}, therefore, implies that consistent estimation of portfolio weights is possible even in the case of growing exposure.

Were we instead to assume that $\bar{s}= O(1)$ so that the number of non-zero elements in each row of the inverse variance matrix was finite, we would have $\|w_u \|_1 = O(1)$. Theorem \ref{thm3} would then become $\|\hat{w}_u - w_u\|_1 = O_p(\sqrt{\mathrm{log}p/n}) = o_p (1)$. This implies that we would obtain a better approximation by assuming finite gross exposure.

\item Using a factor model  to estimate the weight, \cite{fan15} derive the same rate as that referred to  in Remark 3 under the condition that the maximum number of non-zero elements  in each row of $\Sigma$ is finite.
\end{enumerate}

We now turn to the estimation of the weights of the \cite{markowitz52} portfolio. The well-known solution of the Markowitz portfolio optimization problem is
\begin{equation}
w^*=\dfrac{D-\rho_1 B}{AD-B^{2}}(\Sigma^{-1}1_{p}/p)+\dfrac{\rho_1 A-B}{AD-B^{2}}(\Sigma^{-1}\mu/p). \label{marksol}
\end{equation}
Since $ \Sigma $ is positive-definite, $ A>0 $ and $ D>0 $, according to the Cauchy-Schwarz inequality, the system has a solution if $ AD-B^{2}>0 $.

 \noindent The optimal weight vector $w^*$ is estimated by
\begin{equation*}
  \hat{w} = \frac{\hat{D} - \rho_1 \hat{B}}{\hat{A}\hat{D} - \hat{B}^2}(\hat{\Theta} 1_p/p) + \frac{\rho_1  \hat{A} - \hat{B}}{\hat{A} \hat{D} - \hat{B}^2}( \hat{\Theta} \hat{\mu}/p).
\end{equation*}

\noindent Theorem \ref{thm4} establishes the consistency of the estimated weights.
\begin{theorem}\label{thm4}
{\it Under Assumptions 1-3, with $\bar{s}^{3/2} \sqrt{logp/n}= o(1)$, assuming $\lambda_j = O ( \sqrt{\mathrm{log}p/n})$ uniformly in j, with $ (AD-B^2) \ge C_1 >0$,\ where $C_1$ is a positive constant, and with $\rho_1$ being uniformly bounded away from infinity, we have
\[\| \hat{w} - w^* \|_1 = O_p (   \bar{s}^{3/2}  \sqrt{logp/n}) = o_p (1).\]
}
\end{theorem}

{\bf Remarks:}
\begin{enumerate}
\item We can show that
\[ \| w^* \|_1 \le \frac{|D - \rho_1 B| \|\Theta 1_p/p \|_1 }{ C_1}   + \frac{|\rho_1 A - B| \|\Theta \mu/p \|_1}{ C_1} = O (\sqrt{\bar{s}}),    \]
implying that we allow for growing exposure. To establish this result, we use $\|\Theta 1_p/p \|_1 = O(  \sqrt{\bar{s}})$ and $\|\Theta \mu/p \|_1 = O (  \sqrt{\bar{s}})$ from Theorem A.1(iii). To bound the other terms in the numerator, we use Lemma A.4 which shows that $A= O(1), D = O(1), |B|=O(1)$. By assumption, $\rho_1$ is bounded and $AD - B^2  \ge C_1$ which completes the proof.

\item Assuming $\bar{s}= O(1)$ leads to finite gross exposure, and the rate in Theorem \ref{thm4} becomes
\[ \|\hat{w} - w^* \|_1 = O_p ( \sqrt{\mathrm{log}p/n}) = o_p (1).\]
The rate of approximation of the optimal portfolio improves greatly under this stronger condition. This result shows that even for very large portfolios we can accurately estimate the weights.

\item With a non-sparse $\Sigma^{-1}$ and with growing exposure, the rate in Theorem \ref{thm4} becomes
\[ \| \hat{w} - w^* \|_1 = O_p ( p^{3/2} \sqrt{logp/n}),\]
which means that we need $p^{3/2} \sqrt{logp/n}=o(1)$ or, equivalently, $ p \log \left(p^{1/3}\right) = o(n^{1/3})$. With a non-sparse $\Sigma^{-1}$, $p$ has to be much smaller than $n$ to maintain consistency.
\end{enumerate}

\subsection{Portfolio Risk Estimation Error}\label{subsec:risk}

We now turn to the portfolio risk estimation error defined as $|\hat{w}_u' (\hat{\Sigma} - \Sigma) \hat{w}_u|$ and $|\hat{w}' (\hat{\Sigma} - \Sigma) \hat{w}|$  for the global minimum variance portfolio and the Markowitz portfolio, respectively. 

The following theorem shows that the risk estimation error for a large global minimum variance portfolio with growing exposure, converges in probability to zero.

\begin{theorem}\label{rthm1}
 {\it Under Assumptions 1-3, $(\bar{s})^{3/2} \sqrt{\mathrm{log}p/n}=o(1)$, and with $\lambda_j = O (\sqrt{logp/n})$ uniformly in j we have
\[ | \hat{w}_u' (\hat{\Sigma} - \Sigma) \hat{w}_u| = O_p ( \bar{s} \sqrt{\mathrm{log}p/n})=o_p(1).\]}
\end{theorem}

This novel result highlights the trade-off between the number of assets, the sample size, and $\bar{s} $. Assuming constant exposure $\|w_u \|_1 = O (1)$, the rate in Theorem \ref{rthm1} becomes
\[ | \hat{w}_u' (\hat{\Sigma} - \Sigma) \hat{w}_u| = O_p ( \sqrt{\mathrm{log}p/n})=o_p(1),\] which is similar to the result reported by  \cite{fan15} in a factor model context.

Theorem \ref{rthm2} is the counterpart of Theorem \ref{rthm1} for the Markowitz portfolio.

\begin{theorem}\label{rthm2}
{\it Under Assumptions 1-3, with $(\bar{s})^{3/2} \sqrt{logp/n}= o(1)$, assuming $\lambda_j = O ( \sqrt{\mathrm{log}p/n})$ uniformly in j, with $ (AD-B^2) \ge C_1 >0$ where $C_1$ is a positive constant, and with $\rho_1$ uniformly bounded away from infinity we have
\[| \hat{w}' (\hat{\Sigma} - \Sigma) \hat{w}  | = O_p (   \bar{s}  \sqrt{logp/n}) = o_p (1).\]
}
\end{theorem}

This result shows that $\bar{s}$ must be small to establish consistency.

\subsection{Ensuring Positive Definiteness of $\hat \Theta$}

While Lemma A.1 shows that $\hat \Theta$ will be positive definite with high probability, and while we have never encountered zero or negative variances in our simulations and empirical application, it could still occur that $\hat \Theta$ is not positive definite in finite samples. Positive definiteness of $\hat{\Theta}$ can be ensured by symmetrization  as in \cite{fane18}, and then eigenvalue cleaning as in \cite{ckm2016}, and  \cite{hautsch2012blocking}. Note that we do not use this procedure anywhere in our simulations or applications.

\begin{enumerate}
  \item Use the matrix symmetrization procedure of \cite{fane18}: let $\hat{\Theta}_{jk}$ represent $j$th row and $k$ th column of $p \times p$ matrix. Construct the symmetric matrix as
\[ \hat{\Theta}_{jk}^s = \hat{\Theta}_{jk} 1_{ \{ | \hat{\Theta}_{jk}| \le | \hat{\Theta}_{kj} | \} } +
\hat{\Theta}_{kj} 1_{ \{ | \hat{\Theta}_{jk}| > | \hat{\Theta}_{kj} | \} }.\]

  \item Use eigenvalue cleaning as in \cite{ckm2016}, and \cite{hautsch2012blocking} to make $\hat{\Theta}^s$ positive definite. Compute the spectral decomposition of  $\hat{\Theta}^s$. In the diagonal matrix of eigenvalues, replace all the eigenvalues that are below some small positive threshold by the value of this threshold. Use the eigenvector matrix and the new diagonal eigenvalues matrix to reconstruct $\hat{\Theta}^{pd,s}$ from the spectral decomposition formula.

  \item Use $\hat{\Theta}^{pd,s}$, which is symmetric and positive definite, instead of $\hat{\Theta}$ when there are zero or negative variances of the portfolio.
\end{enumerate}

\section{Simulations} \label{sec:simulation}

This section begins with a discussion of  the implementation of the nodewise estimator as well as alternative estimators. We then present the setup of this simulation study after which we report and discuss our results.

\subsection{Implementation of the Nodewise Estimator}\label{GIC}
The nodewise regression approach is implemented by using the coordinate descent algorithm in the \texttt{glmnet} package \citep{glmnet}.  The tuning parameters $\lambda_j$ are chosen using the GIC \citep{fan13} as defined in the subsection  \ref{method}.  One important fact to note is that in our method, we have $p$ tuning parameters which may be larger than the sample size. This is a difficult issue.  However, our estimation of the rows of precision matrix is not affected by this problem, as can be seen in Theorem A.1. Furthermore, our theorem statements in the main text here, use uniformity in $j=1,\cdots, p$ for $\lambda_j$.  Note that when we estimate the rows of the precision matrix, there is only one tuning parameter estimated for each row, with sample size $n$. This process is repeated $p$ times to get the precision matrix estimate. In the case of $n>p$, we do not have any issues. From simulations, and an out-of-sample forecasting exercise, we have not found  any effects of this problem affecting our results.


   
\subsection{Alternative Covariance Matrix Estimation Methods}\label{sec:exist}

\subsubsection{Ledoit-Wolf Shrinkage}

\cite{ledoit04} propose an estimator for high-dimensional covariance matrices that is invertible and well-conditioned. Their estimator is a linear combination of the sample covariance matrix and an identity matrix, and they demonstrate that their estimator is asymptotically optimal for the quadratic loss function.

\subsubsection{Multi-factor Estimator}
The Arbitrage Pricing Theory that is derived by \cite{ross76,ross77}, and the multi-factor models that are proposed by \cite{chamberlain83} have motivated the use of factor models for the estimation of excess return covariance matrices. The model takes the form
\begin{equation}
\label{facmodel}
\mathbf{r}_F=\mathbf{B} \mathbf{f}+\mathbf{\varepsilon},
\end{equation}
where $\mathbf{r}_F$ is a matrix of dimensions $p\times n$ of excess returns of the assets over the risk-free interest rate, and $\mathbf{f}$ is a $K\times n$ matrix of factors. $\mathbf{B}$ is a $p\times K$ matrix of unknown factor loadings, and $\mathbf{\varepsilon}$ is the matrix of idiosyncratic error terms uncorrelated with $\mathbf{f}$. This model yields an estimator for the covariance matrix of ${\bf r}_F$:
\begin{equation}
\Sigma_{FAC}=\mathbf{B} \mathrm{cov}(\mathbf{f})\mathbf{B}^{\prime}
+\Sigma_{n,0},
\end{equation}
where $\Sigma_{n,0}$ is the covariance matrix of errors $\varepsilon$.
When the factors are observed, as in \cite{fan08}, the matrix of loadings $\mathbf{B}$ can be estimated by least squares. In the case where the factors are unobserved, \cite{fan2013large,fan15} proposed the POET to estimate $\Sigma_{FAC}$ and then invert that to get the precision matrix estimate.

\subsubsection{Positive Semi-definite Projection}

This method by \cite{zhao14} begins with Kendall's tau estimator, instead of the sample covariance matrix. It uses the following constrained optimization to get a positive semi-definite version of Kendall's tau estimator:
\[ \tilde{S} = argmin_S \| \hat{S}_{kt} - S \|_{\infty} \quad {\mbox s.to}  \quad S \ge 0,\]
where $\hat{S}_{kt}$ is Kendall's tau estimator defined in section 2.1 of \cite{zhao14}. 
In the simulations that follow, we invert $\tilde{S}$ and use it in our global minimum variance and the Markowitz portfolios. This method requires the setting of many tuning parameters and is very sensitive to their choice. We used the values suggested by \cite{zhao14}. In limited simulations that we have done, the estimator of variance of the portfolio under this method varies widely, and  does not  converge to its theoretical value. Hence, we do not include this method in our Tables.

\subsection{Simulation Setup and Results}

\subsubsection{DGP}

For our simulations, we use two pairs of data generating processes (DGPs). In the first pair, the covariance matrix of the returns has a Toeplitz structure. For each $t=1,2\cdots, n$, the $p\times 1$ vector of excess returns  $r_t$ follows multivariate normal distribution with covariance matrix $\Omega$. The entry on row $k\in 1,\cdots, p$, and column $l\in 1,\cdots p$ of $\Omega$ has value $\omega^{|k-l|}$. We set $\omega=0.15$. Setting $\omega=0.25$ yielded similar results which we therefore, do not report. In the first DGP of this pair, the mean vector of the excess returns $\mu=0_p$. We use this DGP with the Global Minimum Variance portfolio and report the results in Tables 1a-b.

In the second DGP of this pair, the mean vector $\mu$ is drawn from a normal distribution mean of zero and variance 0.0001. This DGP is used with the \cite{markowitz52} portfolio. The return target is set at $\rho_1=0.000376$, corresponding to a 10\% annual return over 252 trading days. The results are reported in Tables 2a-2b.

In the second pair of DGPs, the excess returns follow a weakly sparse 3 factor structure. In the first DGP of this pair, and for each $t=1,\cdots, n$
\[ r_t= B_t f_t + e_t,\]
where $e_t\sim \mathcal{N}(0,I_p)$.The vector of factors $f_t$ is drawn from $\mathcal{N}(0_3, I_3/10)$. The entries of the matrix of factor loadings for asset $j=1,\cdots, p$, $B_{t,j}$, are drawn from $\mathcal{N}(0_3, I_3/100)$, and the full matrix of loadings is given by $B_t= (B_{t,1}, B_{t,2}, \cdots B_{t,j}, \cdots, B_{t,p})'$.  We use this DGP with the Global Minimum Variance portfolio and report the results in Tables 3a-b.

The second DGP of this pair adds a non-zero mean to the excess returns. For $t=1,\cdots, n$
\[ r_t = \mu+ B_t f_t + e_t,\]
where $\mu\sim \mathcal{N}(0, 0.01)$. We use this DGP with the \cite{markowitz52} portfolio and report the results in Tables 4a-4b.


\subsubsection{Dimensions}

For the four DGPs and the four estimators described above,  we simulate data with both $p=n/2$ and $p=3*n/2$. In both cases, the number of assets grows with the sample size, but in one case $p<n$, and in the other $p>n$. We use sample sizes $n=100, 200, 400$. Due to the computation cost of evaluating the four methods on four DGPs, we limit ourselves to 100 replications.

\subsubsection{Metrics and Methods}

For each experiment, we consider three metrics: the variance estimation error, the weight estimation error, and the risk estimation error. These metrics are headers, "Variance", "Weight", "Risk" in the Tables.
These metrics are used in Theorems 3.1, 3.3, and 3.5 (respectively) in the case of the Global Minimum Variance portfolios, and in Theorems 3.2, 3.4, and 3.6 in the case of the Markowitz portfolios. Nodewise, POET, Ledoit-Wolf, in the first column of each table, refer to methods that are analyzed. Nodewise is our method in the paper; POET refers to factor-model based estimation in  subsection 4.2.2; Ledoit-Wolf refers to the shrinkage method in subsection 4.2.1.

\subsubsection{Performance}

In all tables the variance, weight, and risk metrics for the nodewise estimator converge to zero when the sample size increases, as predicted by the theory. For instance, in Table 1a the risk estimation error goes from $0.0038$ with $n=100, p=50$ to $0.0003$ with $n=400, p=200$, and in Table 2a, the variance estimation error goes from $0.3571$ with $n=100, p=50$ to $0.2750$ with $n=400, p=600$.

 Compared to other methods, Nodewise provides the best result in terms of risk estimation error in Tables 1a-b, 2a-2b, with Toeplitz  DGP. For example with $n=200, p=300$ in Table 2b, Nodewise has a risk error of 0.0013 which is lower than Ledoit-Wolf at 0.0018, POET at 0.0019. For variance and weight estimation errors, Nodewise and Ledoit-Wolf methods perform well compared to other two methods.
 


\begin{table}[!htp]
\begin{centering}
\resizebox{\textwidth}{!} {
\begin{tabular}{l ccc ccc ccc }

 \multicolumn{10}{c}{Table 1a: TOEPLITZ DGP: GLOBAL MINIMUM VARIANCE PORTFOLIO, $p=n/2$}\\
\midrule
& \multicolumn{3}{c}{$n=100$, $p=50$} & \multicolumn{3}{c}{$n=200$, $p=100$} & \multicolumn{3}{c}{$n=400$, $p=200$} \\

 & Variance & Weight & Risk & Variance & Weight & Risk & Variance & Weight & Risk \\
\cmidrule(lr){2-4} \cmidrule(lr){5-7} \cmidrule(lr){8-10}
  Nodewise & 0.4013& 0.2488 & 0.0038 & 0.3788 & 0.1718 & 0.0012 & 0.3624&0.1180&0.0003 \\ 
  POET & 0.5695 & 0.3874 & 0.0146 & 0.4222 & 0.2869 & 0.0054 & 0.3281&0.2061&0.0029 \\ 
  Ledoit-Wolf & 0.3216 & 0.0642 & 0.0066 & 0.3173 & 0.0602 & 0.0033 & 0.3200 & 0.0572&0.0017 \\ 
 \bottomrule
\end{tabular}
}
\end{centering}
\end{table}

\begin{table}[!htp]
\begin{centering}
\resizebox{\textwidth}{!} {
\begin{tabular}{l ccc ccc ccc }

 \multicolumn{10}{c}{Table 1b: TOEPLITZ DGP: GLOBAL MINIMUM VARIANCE PORTFOLIO, $p=3*n/2$}\\
\midrule
& \multicolumn{3}{c}{$n=100$, $p=150$} & \multicolumn{3}{c}{$n=200$, $p=300$} & \multicolumn{3}{c}{$n=400$, $p=600$} \\

 & Variance & Weight & Risk & Variance & Weight & Risk & Variance & Weight & Risk \\
\cmidrule(lr){2-4} \cmidrule(lr){5-7} \cmidrule(lr){8-10}
  Nodewise & 0.4185& 0.2339 & 0.0013 & 0.3883 & 0.1628 & 0.0004 & 0.3697&0.1155&0.0001 \\ 
  POET & 0.6616 & 0.3400 & 0.0048 & 0.4617 & 0.2653 & 0.0018 & 0.3237&0.2252&0.0007 \\ 
  Ledoit-Wolf & 0.3492 & 0.0516 & 0.0023 & 0.3429 & 0.0421 & 0.0011 & 0.3421 & 0.0373&0.0006 \\ 
 \bottomrule
\end{tabular}
}
\end{centering}
\end{table}

\begin{table}[!htp]
\begin{centering}
\resizebox{\textwidth}{!} {
\begin{tabular}{l ccc ccc ccc }

 \multicolumn{10}{c}{Table 2a: TOEPLITZ DGP: MARKOWITZ PORTFOLIO, $p=n/2$}\\
\midrule
& \multicolumn{3}{c}{$n=100$, $p=50$} & \multicolumn{3}{c}{$n=200$, $p=100$} & \multicolumn{3}{c}{$n=400$, $p=200$} \\

 & Variance & Weight & Risk & Variance & Weight & Risk & Variance & Weight & Risk \\
\cmidrule(lr){2-4} \cmidrule(lr){5-7} \cmidrule(lr){8-10}
  Nodewise & 0.3571& 0.3920 & 0.0118 & 0.2818 & 0.3031 & 0.0036 & 0.2750 &0.1875&0.0017 \\ 
  POET & 0.5074 & 0.4862 & 0.0148 & 0.3744 & 0.2688 & 0.0056 & 0.3298&0.2088&0.0024 \\ 
  Ledoit-Wolf & 0.3416 & 0.3921 & 0.0071 & 0.2722 & 0.1759 & 0.0035 & 0.2704 & 0.1494&0.0017 \\ 
 \bottomrule
\end{tabular}
}
\end{centering}
\end{table}

\begin{table}[!htp]
\begin{centering}
\resizebox{\textwidth}{!} {
\begin{tabular}{l ccc ccc ccc }

 \multicolumn{10}{c}{Table 2b: TOEPLITZ DGP: MARKOWITZ PORTFOLIO, $p=3*n/2$}\\
\midrule
& \multicolumn{3}{c}{$n=100$, $p=150$} & \multicolumn{3}{c}{$n=200$, $p=300$} & \multicolumn{3}{c}{$n=400$, $p=600$} \\

 & Variance & Weight & Risk & Variance & Weight & Risk & Variance & Weight & Risk \\
\cmidrule(lr){2-4} \cmidrule(lr){5-7} \cmidrule(lr){8-10}
  Nodewise & 0.2975& 0.2836 & 0.0024 & 0.2825 & 0.1821 & 0.0013 & 0.2696&0.1086&0.0006 \\ 
  POET & 0.4405 & 0.3198 & 0.0048 & 0.3727 & 0.2580 & 0.0019 & 0.3289&0.2001&0.0008 \\ 
  Ledoit-Wolf & 0.2693 & 0.1288 & 0.0024 & 0.2674 & 0.1020 & 0.0018 & 0.2616 & 0.0425&0.0006 \\ 
 \bottomrule
\end{tabular}
}
\end{centering}
\end{table}

\begin{table}[!htp]
\begin{centering}
\resizebox{\textwidth}{!} {
\begin{tabular}{l ccc ccc ccc }

 \multicolumn{10}{c}{Table 3a: SPARSE FACTOR DGP: GLOBAL MINIMUM VARIANCE PORTFOLIO, $p=n/2$}\\
\midrule
& \multicolumn{3}{c}{$n=100$, $p=50$} & \multicolumn{3}{c}{$n=200$, $p=100$} & \multicolumn{3}{c}{$n=400$, $p=200$} \\

 & Variance & Weight & Risk & Variance & Weight & Risk & Variance & Weight & Risk \\
\cmidrule(lr){2-4} \cmidrule(lr){5-7} \cmidrule(lr){8-10}
  Nodewise & 0.0760& 0.2226 & 0.0027 & 0.0415 & 0.1605 & 0.0009 & 0.0179&0.1129&0.0003 \\ 
  POET & 0.3593 & 0.2886 & 0.0080 & 0.1718 & 0.2626 & 0.0021 & 0.0895&0.1477&0.0006 \\ 
  Ledoit-Wolf & 0.0171 & 0.0188 & 0.0003 & 0.0090 & 0.0155 & 0.0001 & 0.0042 & 0.0163&0.0001 \\ 
 \bottomrule
\end{tabular}
}
\end{centering}
\end{table}

\begin{table}[!htp]
\begin{centering}
\resizebox{\textwidth}{!} {
\begin{tabular}{l ccc ccc ccc }

 \multicolumn{10}{c}{Table 3b: SPARSE FACTOR DGP: GLOBAL MINIMUM VARIANCE PORTFOLIO, $p=3*n/2$}\\
\midrule
& \multicolumn{3}{c}{$n=100$, $p=150$} & \multicolumn{3}{c}{$n=200$, $p=300$} & \multicolumn{3}{c}{$n=400$, $p=600$} \\

 & Variance & Weight & Risk & Variance & Weight & Risk & Variance & Weight & Risk \\
\cmidrule(lr){2-4} \cmidrule(lr){5-7} \cmidrule(lr){8-10}
  Nodewise & 0.0800& 0.2285 & 0.0010 & 0.0373 & 0.1617 & 0.0003 & 0.0159&0.1148&0.00009 \\ 
  POET & 0.3221 & 0.2921 & 0.0025 & 0.1762 & 0.2331 & 0.0008 & 0.1008&0.1960&0.00030 \\ 
  Ledoit-Wolf & 0.0138 & 0.0225 & 0.0001 & 0.0065 & 0.0209 & 0.00002 & 0.0027 & 0.0198&0.00001 \\ 
 \bottomrule
\end{tabular}
}
\end{centering}
\end{table}

\begin{table}[!htp]
\begin{centering}
\resizebox{\textwidth}{!} {
\begin{tabular}{l ccc ccc ccc }

 \multicolumn{10}{c}{Table 4a: SPARSE FACTOR DGP: MARKOWITZ PORTFOLIO, $p=n/2$}\\
\midrule
& \multicolumn{3}{c}{$n=100$, $p=50$} & \multicolumn{3}{c}{$n=200$, $p=100$} & \multicolumn{3}{c}{$n=400$, $p=200$} \\

 & Variance & Weight & Risk & Variance & Weight & Risk & Variance & Weight & Risk \\
\cmidrule(lr){2-4} \cmidrule(lr){5-7} \cmidrule(lr){8-10}
  Nodewise & 0.1714& 0.4596 & 0.0027 & 0.0367 & 0.2425 & 0.0009 & 0.0279 &0.1858&0.00030 \\ 
  POET & 0.3066 & 0.3928 & 0.0085 & 0.1509 & 0.2210 & 0.0021 & 0.0821&0.1573&0.00060 \\ 
  Ledoit-Wolf & 0.0866 & 0.2592 & 0.0004 & 0.0166 & 0.0890 & 0.0001 & 0.0081 & 0.0551&0.00002 \\ 
 \bottomrule
\end{tabular}
}
\end{centering}
\end{table}

\begin{table}[!htp]
\begin{centering}
\resizebox{\textwidth}{!} {
\begin{tabular}{l ccc ccc ccc }

 \multicolumn{10}{c}{Table 4b: SPARSE FACTOR DGP: MARKOWITZ PORTFOLIO, $p=3*n/2$}\\
\midrule
& \multicolumn{3}{c}{$n=100$, $p=150$} & \multicolumn{3}{c}{$n=200$, $p=300$} & \multicolumn{3}{c}{$n=400$, $p=600$} \\

 & Variance & Weight & Risk & Variance & Weight & Risk & Variance & Weight & Risk \\
\cmidrule(lr){2-4} \cmidrule(lr){5-7} \cmidrule(lr){8-10}
  Nodewise & 0.0586& 0.2618 & 0.0009 & 0.0294 & 0.1901 & 0.0003 & 0.0098&0.1234&0.00009 \\ 
  POET & 0.2396 & 0.3029 & 0.0025 & 0.1473 & 0.2354 & 0.0008 & 0.0925&0.1961&0.00030 \\ 
  Ledoit-Wolf & 0.0157 & 0.0739 & 0.0001 & 0.0075 & 0.0484 & 0.00002 & 0.0028 & 0.0246&0.00001 \\ 
 \bottomrule
\end{tabular}
}
\end{centering}
\end{table}

\section{Empirical application}\label{sec:per}

\subsection{Performance Measures}

 In this section, we compare the Nodewise regression approach to the POET estimator and the \cite{ledoit04} estimator in an out-of-sample portfolio optimization application. We focus on four metrics commonly used in finance, namely, the Sharpe ratio (SR from now on), the portfolio turnover, and the average return and variance of portfolios. We consider portfolio formation with and without transaction costs.


We use a rolling horizon method for out-of-sample forecasting. Samples of length $n$ are split into an in-sample training part indexed $(1:n_I)$, and an out-of-sample testing part indexed $(n_I +1:n)$. The rolling window method works as follows: the portfolio weights $\hat{w}_{n_I}$ are calculated in-sample for the period in between $(1: n_I)$, then multiplied by the return at time $n_I+1$ to get the out-of-sample portfolio returns $\hat{w}_{n_I}' r_{n_I +1}$. We roll the window by one period $(2: n_{I}+1)$ and form the portfolio weight $\hat{w}_{n_I +1}$ for that period. This weight vector is multiplied by the returns at time $n_I +2$ to get the portfolio returns $\hat{w}_{n_I+1}' r_{n_I +2}$. We repeat this procedure until the end of the sample.

With no transaction costs, the out-of-sample average portfolio return and variance are
\begin{align*}
\hat{\mu}_{os} = \frac{1}{n-n_I} \sum_{t=n_I}^{n-1} \hat{w}_t' r_{t+1} & \text{ and  } \hat{\sigma}^2_{os} = \frac{1}{(n-n_I)-1} \sum_{t=n_I}^{n-1} (\hat{w}_t' r_{t+1} - \hat{\mu}_{os})^2.
\end{align*}
Using these two statistics we compute the Sharpe ratio:
\[ SR = \hat{\mu}_{os}/\hat{\sigma}_{os}.\]

When transaction costs are present, the definition of the metrics is modified. Let $c$ be the transaction cost, which we set to 50 basis points throughout following \cite{demiguel09b}. From \cite{ban16, li15}, the excess portfolio return at time $ t $ with transaction cost is
\[ Return_t = \hat{w}_t' r_{t+1} - c ( 1 + \hat{w}_t' r_{t+1}) \sum_{j=1}^p | \hat{w}_{t+1,j} - \hat{w}_{t,j}^+ |,\]
where $\hat{w}_{t,j}^+ =  \hat{w}_{t,j}  (1+ R_{t+1,j})/(1+R_{t+1,p})$, and $R_{t+1,j}$ is the excess return added to the risk-free rate for $j^{th}$ asset, and
$R_{t+1,p}$ is the portfolio excess return plus risk-free rate.

Using this definition of the excess portfolio returns, we define the mean and variance of the portfolio as
\begin{align*}
\hat{\mu}_{os,c} = \frac{1}{n-n_I} \sum_{t=n_I}^{n-1}  Return_t & \text{ and  } \hat{\sigma}^2_{os,c} = \frac{1}{(n-n_I)-1} \sum_{t=n_I}^{n-1} (Return_t  - \hat{\mu}_{os,c})^2,
\end{align*}
and the SR with transaction costs becomes:
\[ SR_c =    \hat{\mu}_{os,c}/ \hat{\sigma}_{os,c}.\]
Finally, we also consider the portfolio turnover:
 \[ PT= \frac{1}{n- n_I} \sum_{t= n_I}^{n-1} \sum_{j=1}^p | \hat{w}_{t+1,j} - \hat{w}_{t,j}^+ |.\]

\subsection{Data}\label{ssbenc}

We use daily and monthly returns of components of the S\&P500 index, and take the 3-month Treasury bill rate as our measure of the risk-free return. With the monthly data we have $n<p$ while for daily data we have $ p<n $. We use two different out-of-sample periods each for daily and for monthly data. We have also tried other possibilities and sub-intervals as reported in the Appendix. The results are similar.

\begin{enumerate}

\item{\bf Monthly Data:} January 1994 to May 2018 with $n=293 $ and $ p=304 $.

 a) In-Sample 1: January 1994-March 2010 ($n_I=195$), Out-Of-Sample 1: April 2010-May 2018 ($n-n_I=98$).

 b) In-Sample 2: January 1994-May 2008 ($n_I=173$), Out-Of-Sample 2: June 2008-May 2018 ($n-n_I= 120$).

\item {\bf Daily Data:} July 2, 2013 to April 30, 2018 with $n=1216 $ and $ p=452 $.

a) In-Sample 1: July 2, 2013-April 28, 2017 ($n_I=964$), Out-Of-Sample 1: May 1, 2017-April 30, 2018 ($n-n_I= 252$).

b) In-Sample 2: July 2, 2013-January 29, 2018 ($n_I=1153$), Out-Of-Sample 2: January 30, 2018-April 30, 2018 ($n-n_I=63$).

\end{enumerate}

For each new training window, the portfolios are rebalanced, and the expected return vector and the covariance matrices are re-estimated. For instance, for a ten-year ($n-n_I=120$) rolling window forecast horizon, we estimate expected returns and covariance matrices, and formulate the Global Minimum  Variance and the Markowitz portfolio $ 120 $ times. Portfolios are held for one month and rebalanced at the beginning of the next month. As a return target in the Markowitz portfolio, we use a monthly target of 0.7974\% and a daily target of 0.0378\%  both of which are equivalents of 10\% yearly return when compounded.

\subsection{Results}
We report the Global Minimum  Variance and the Markowitz portfolio empirical results with and without transaction costs (TC in tables) based on the POET, Nodewise and Ledoit-Wolf estimators.

In Table 5, we  report the monthly portfolio performances for the full evaluation period of January 1994 - May 2018.
 Nodewise regression-based portfolios have the highest SR and the lowest variance and turnover rates in all cases. Portfolios based on the Ledoit-Wolf estimator generally yield the highest returns, but at the price of a ten-fold increase in variance leading to the smallest SR.

\begin{table}
\begin{centering}
\resizebox{\textwidth}{!} {
\begin{tabular}{rrrcccrrcc}\\
\multicolumn{10}{c}{\large \bf Table 5: Monthly Portfolio Returns-Variance-Sharpe Ratio-Turnover}\\
\toprule
\multicolumn{1}{c}{}&\multicolumn{4}{c}{\bfseries Global Minimum Portfolio}&\multicolumn{1}{c}{\bfseries }&\multicolumn{4}{c}{\bfseries Markowitz Portfolio}\tabularnewline
\cline{2-5} \cline{7-10}
\multicolumn{1}{r}{}&\multicolumn{1}{c}{Return}&\multicolumn{1}{c}{Variance}&\multicolumn{1}{c}{SR}&\multicolumn{1}{c}{Turnover}&\multicolumn{1}{c}{}&\multicolumn{1}{c}{Return}&\multicolumn{1}{c}{Variance}&\multicolumn{1}{c}{SR}&\multicolumn{1}{c}{Turnover}\tabularnewline
\midrule
\multicolumn{10}{c}{\bf In-Sample: January 1994-March 2010, Out-Of-Sample: April 2010-May 2018, $ n_I=195 $, $ n-n_I=98 $}\tabularnewline
{\bf without TC }&&&&&&&&&\tabularnewline
POET&0.02499&0.01953&0.1788&0.2080&& 0.03317&0.05248& 0.1448&-\tabularnewline
Nodewise&0.02644&0.01580&0.2104&0.1268&& 0.03503&0.04585& 0.1635&-\tabularnewline
Ledoit-Wolf&0.07140&0.18421&0.1664&0.2150&&-0.01128&0.24568&-0.0227&-\tabularnewline
\midrule
{\bf with TC}&&&&&&&&&\tabularnewline
POET&0.02497&0.01964&0.1782&    -&& 0.03324&0.05274& 0.1447&   0.2497\tabularnewline
Nodewise&0.02650&0.01590&0.2102&    -&& 0.03517&0.04613& 0.1637&  0.1709  \tabularnewline
Ledoit-Wolf&0.07175&0.18545&0.1666&    -&&-0.01166&0.24733&-0.0234&0.4421\tabularnewline
\midrule
\hline
\multicolumn{10}{c}{\bf In-Sample: January 1994-May 2008, Out-Of-Sample: June 2008-May 2018 $ n_I=173 $, $ n-n_I=120 $ }\tabularnewline
{\bf without TC}&&&&&&&&&\tabularnewline
POET&0.02078&0.01805&0.1547&0.1806&& 0.022689&0.02537& 0.1424&-\tabularnewline
Nodewise&0.02158&0.01254&0.1927&0.1357&& 0.024571&0.02054& 0.1714&-\tabularnewline
Ledoit-Wolf&0.05039&0.11844&0.1464&0.2795&&-0.009022&0.13209&-0.0248&-\tabularnewline
\midrule
{\bf with TC}&&&&&&&&&\tabularnewline
POET&0.02070&0.01813&0.1538&    -&& 0.022634&0.02548& 0.1418&   0.2232 \tabularnewline
Nodewise&0.02156&0.01261&0.1920&    -&& 0.024580&0.02064& 0.1710&  0.1759  \tabularnewline
Ledoit-Wolf&0.05041&0.11904&0.1461&    -&&-0.009442&0.13282&-0.0259&  1.8899  \tabularnewline
\bottomrule \hline
\end{tabular}
} 
\end{centering}
\end{table}

\begin{table}

\begin{centering}

\resizebox{\textwidth}{!} {
\begin{tabular}{rrrcccrrcc}\\
\multicolumn{10}{c}{\large \bf Table 6: Daily Portfolio Returns-Variance-Sharpe Ratio-Turnover}\\
\toprule
\multicolumn{1}{c}{}&\multicolumn{4}{c}{\bfseries Global Minimum Portfolio}&\multicolumn{1}{c}{\bfseries }&\multicolumn{4}{c}{\bfseries Markowitz Portfolio}\tabularnewline
\cline{2-5} \cline{7-10}
\multicolumn{1}{r}{}&\multicolumn{1}{c}{Return}&\multicolumn{1}{c}{Variance}&\multicolumn{1}{c}{SR}&\multicolumn{1}{c}{Turnover}&\multicolumn{1}{c}{}&\multicolumn{1}{c}{Return}&\multicolumn{1}{c}{Variance}&\multicolumn{1}{c}{SR}&\multicolumn{1}{c}{Turnover}\tabularnewline
\midrule
\multicolumn{10}{c}{\bf In-Sample: Jul 2 2013-Apr 28 2017, Out-Of-Sample: May 1 2017-Apr 30 2018, $ n_I=964 $, $ n-n_I=252 $ }\tabularnewline
{\bf without TC}&&&&&&&&&\tabularnewline
POET& 4.537e-04&4.757e-05& 0.0657&0.0685&&3.757e-04&4.630e-05&0.0552&-\tabularnewline
Nodewise& 4.424e-04&4.499e-05& 0.0659&0.0570&&3.731e-04&4.421e-05&0.0561&-\tabularnewline
Ledoit-Wolf& 4.103e-04&3.566e-05& 0.0687&0.3430&&5.329e-04&3.738e-05&0.0871&-\tabularnewline
\midrule
{\bf with TC}&&&&&&&&&\tabularnewline
POET& 4.225e-04&4.765e-05& 0.0612& -&&3.268e-04&4.629e-05&0.0480&0.0884\tabularnewline
Nodewise& 4.161e-04&4.492e-05& 0.0620& -&&3.319e-04&4.416e-05&0.0499&0.0721\tabularnewline
Ledoit-Wolf& 4.174e-05&3.563e-05& 0.0069& -&&1.490e-04&3.766e-05&0.0242&0.3539\tabularnewline
\midrule \hline
\multicolumn{10}{c}{\bf In-Sample: Jul 2 2013-Jan 29 2018, Out-Of-Sample: Jan 30 2018-Apr 30 2018, $ n_I=1153 $, $ n-n_I=63 $ }\tabularnewline
{\bf without TC}&&&&&&&&&\tabularnewline
POET&-9.237e-04&1.330e-04&-0.0800&0.0110&&-0.001035&1.217e-04&-0.0938&-\tabularnewline
Nodewise&-8.794e-04&1.262e-04&-0.0782&0.0540&&-0.000934&1.187e-04&-0.0857&-\tabularnewline
Ledoit-Wolf&-8.008e-04&5.874e-05&-0.1044&0.3269&&-0.000718&6.119e-05&-0.0918&-\tabularnewline
\midrule
{\bf with TC}&&&&&&&&&\tabularnewline
POET&-8.059e-04&1.341e-04&-0.0695&-&&-0.000932&1.226e-04&-0.0842&0.0371\tabularnewline
Nodewise&-8.323e-04&1.276e-04&-0.0736&-&&-0.000901&1.200e-04&-0.0823&0.0709\tabularnewline
Ledoit-Wolf&-1.258e-03&5.942e-05&-0.1632&-&&-0.001210&6.221e-05&-0.1534&0.3419\tabularnewline
\bottomrule
\end{tabular}
} 
\end{centering}
  
\end{table}

In Table 6, we  report the daily portfolio performances. When there are no transaction costs, all three estimators perform quite similarly in the first sample, except for much higher turnover for the Ledoit-Wolf estimator. When transaction costs are taken into account, the SR of Ledoit-Wolf portfolio decreases considerably. This deterioration is driven by the high turnover rates and the lower mean returns. Nodewise regression-based portfolios also perform well here. In the second sample, all estimators have negative returns. The Ledoit-Wolf estimator yields the most negative returns with the smallest variance when transaction costs are taken into account.

In summary, the Nodewise estimator shows the best performance in terms of out-of-sample SR, turnover rate, and variance on monthly data with $ n<p $. As to daily data, the results are more mixed, but the Nodewise estimator generally yields the lowest turnover and highest SR. On daily data, the Ledoit-Wolf estimator generally yields portfolios with a low variance, but also low returns. The turnover of portfolio based on Ledoit-Wolf is generally very high.

\section{Conclusion}\label{sec:conclusion}

In this paper, we analyze the variance, weights, and risk of  a large portfolio when the portfolio is constructed using nodewise regression to estimate the inverse covariance matrix of excess returns of assets. We show that all three quantities can be consistently estimated even when the data are weakly correlated.  We compare our  estimator to the factor model-based and shrinkage estimators in simulations and in an empirical application. Nodewise regression, being based on an exact inverse formula, provides low variance, and good SR  compared with factor model and shrinkage based approaches.

\bibliographystyle{chicagoa}
\bibliography{Nodewise}
\end{document}